\documentclass{amsart}
\usepackage{amsfonts,amssymb,amsmath,amsthm}
\usepackage{url}
\usepackage{enumerate}

\newtheorem{thrm}{Theorem}[section]
\newtheorem{lem}[thrm]{Lemma}
\newtheorem{prop}[thrm]{Proposition}

\theoremstyle{definition}

\numberwithin{equation}{section}

\author{Jelena Stojanov}
\address{
Technical Faculty "Mihajlo Pupin"\\
University of Novi Sad\\
Serbia}
\email{jelena@tfzr.uns.ac.rs}

\keywords{Finsler space, hypersurface, indicatrix, mean curvature}
\subjclass{53B40, 53A10, 14J70,}
\begin{document}

\title{Mean curvature of the indicatrix of Finsler manifold}

\begin{abstract}
Fundamental function in Finsler manifold defines a metrices that depend on a point and a direction. At any point tangent space  is a Riemannian and an indicatrix is a convex hypersurface. In this paper a mean curvature of the indicatrix is expressed in terms of fundamental function.
\end{abstract}
\maketitle

\section{Introduction}
Finsler manifold is generalization of the Riemannian one, in the same way as Riemann manifold is for the Euclidean. A metric depends on the point and the direction.

Let $M$ be an $n-$dimensional differential manifold and $F$ smooth, real nonnegative and $1-$homogeneous function that acts on the tangent bundle of $M$, such that a Hessian of $F^2$ is nonnegative \cite{Rund1959}\cite{Matsumoto1986}\cite{Bao2000}. $F$ is called a fundamental function of Finsler manifold $(M,F)$.

For a fixed point $x\in M$, the tangent vector space $T_xM$ can be considered as a Riemannian space with metric generated by $F$,
$$
g_{ij}=g_{ij}(y)=\dfrac{1}{2}\dfrac{\partial^2 F^2(x,y)}{\partial y^i \partial y^j};
$$
but also,  $T_x M$ is Minkowski space with a norm $F(x,y)$. Collection of  vectors
 $$I_x=\{y \in T_xM | F(x,y)=1 \}$$
 is called an indicatrix, and might be considered twofold: as a hypersurface and as an unit sphere. The first viewpoint appears in \cite{Srivastava1975}\cite{Watanabe1973}\cite{Nishimura1991} and the second one in \cite{Ji2002}\cite{Kikuchi1962}. 
 
 As a Riemannian hypersurface, the indicatrix is convex and orientable \cite{Rund1959}\cite{Matsumoto1986}, and it has the tangent bundle containing  Euclid spaces as fibres. From the theory of hypersurfaces  \cite{Yano1970}, it stems that the indicatrix $I_x$ is totally umbilical with constant unit mean curvature \cite{Kikuchi1962} \cite{Matsumoto1986}.
 
In \cite{Srivastava1975}\cite{Watanabe1973} the indicatrix is represented in certain local frame. Covariant differential of that representation is given, and its integrability conditions. Necessary and sufficient conditions, for any hypersurfaces to be homothetic with the indicatrix, contain a mean curvature as a characteristic of the shape. 

A mean curvature of the indicatrix $I_x$ can be expressed as a trace of second fundamental tensor, arithmetic mean of principal curvatures, but in this paper, it will be connected directly with Hessian of fundamental function. The idea comes from a paper of Nishimura and Hashiguchi \cite{Nishimura1991}, where the Gauss curvature of hypersurface is expressed in terms of its defining function.

\section{Mean curvature of a hypersurface }

Let $I$ be an oriented hypersurface in a Riemannian space $V^n$,
given implicitly by a $C^{\infty}-$differentiable function
$f$,
\begin{equation}\label{def I}
I=\left\{x\in V^n\vert f(x)=0\right\}.
\end{equation}
As usual, notation for coordinate map in Riemannian space is $x=(x^1,x^2,...,x^n)$ and an abbreviation for partial derivatives is used,
$f_i=\partial _i f=\frac{\partial f}{\partial x^i}.$

Differentiability of at least second order provides nonvanishing gradient 
$$\nabla f:I\rightarrow \mathbb{R}^n, \quad \nabla f (x)\neq 0,\qquad
\nabla f=\left[ \nabla f\right]_{n\times 1} =(f_1,f_2,...,f_n).$$

Tangent space on the hypersurface at each point $p\in I$ is a
real vector Euclidean space $T_{p}I$ and has
dimension $n-1$ \cite{Rashevskii1964}\cite{Gudmundsson2010}. So, it is a hyperplane and with the hypersurface has common unit normal vector at the point $p$.
Unit normal vector field on the hypersurface is continuous and globally defined
\begin{equation}\label{normal}
N=\displaystyle\frac{\varepsilon}{\vert \nabla f \vert}\nabla f; \qquad (\varepsilon=\pm1)
\end{equation}

Element of the tangent hyperplane $T_{p}I$ is a vector tangent on
$I$, but also on $V^n$, $v\in T_{p}I\subset T_{p}V^n$, so it has
$n$ components $v=(v^1,v^2,\dots,v^n )$; a basis in $T_{p}I$ has
$n-1$ elements  $\lbrace X_1,X_2,\dots,X_{n-1} \rbrace$ and there
is a decomposition $v=v^{\alpha}X_{\alpha}$. Precisely, Latin indices run $\{1,2,...,n\}$ and mean
components in  $T_{p}V^n$ and Greek ones run $\{1,2,...,n-1\}$ and
consider decomposition over the basis of $T_{p}I$.

Local bending (shape) of the hypersurface is expressed by
infinitesimal changes of the hypersurfaces normal. A shape
operator applied to the tangent space $T_{p}I$ is a linear
transformation and for a tangent vector $v\in T_{p}I$ determines
negative derivative of the normal in the direction of $v$
\cite{Gray2006},
\begin{equation}\nonumber
S:T_{p}I\rightarrow T_{p}I,\qquad S(v)=-D_{v}N.
\end{equation}
Calculation of a directional derivative $D_{v}N$ means partial
derivatives of the normal. They are tangent on the hypersurface
and Wiengartens equations
$\partial_{\beta}N=-h^{\alpha}_{\beta}X_{\alpha}$ allowed matrix
notation
\begin{equation}\nonumber
S(v)=\left[h^{\alpha}_{\beta}\right]_{(n-1)\times(n-1)}\left[v^{\alpha}\right]_{(n-1)\times1}.
\end{equation}

A mean curvature of the hypersurface is arithmetical mean for
diagonal elements in  matrix $\left[h^{\alpha}_{\beta}\right]$ of
normal's derivatives decompositions over the tangent space, i.e. is
proportional to the trace of the shape operator
$$H=\frac{1}{n-1}tr\left[h^{\alpha}_{\beta}\right]=\frac{1}{n-1}trS.$$

$T_{p}I$ is $(n-1)-$dimensional Euclidean space and scalar product
 is defined as an intrinsic form. But,
if one of the arguments is obtained by the shape operator, scalar
product is quadratic form in ambient space $T_{p}V^n$, \cite{Nishimura1991}.
\begin{prop}\label{propozicija}
The shape operator of an oriented hypersurface $(I,N)$, given by
(\ref{def I}), satisfies the following condition
\begin{equation}\nonumber
u\cdot S(v)=-\frac{\varepsilon}{\vert \nabla f
\vert}f_{ij}u^i v^j
\end{equation}
for any vectors $u=(u^i)$ and $v=(v^i)$ tangent on the hypersurface
$I$. An abbreviation is used, $f_{ij}=\partial_{i}\partial_{j}f.$
\end{prop}
\begin{proof}
$\nabla f$ is vector field on an open neighbourhood of $I$, so
directional derivation of (\ref{normal}) gives
$$\partial_{j}(\nabla f)v^j=\varepsilon D_{v}(\vert\nabla f\vert)N-\varepsilon\vert\nabla f\vert S(v).$$
$\partial_{j}(\nabla f)$ is a quadratic matrix of order $n$ with
second partial derivatives. Observing $N,\;S(v)$ and $u$ as
elements of $T_{p}V^n$, scalar product is
$$ u\cdot f_{ij}v^j=-\varepsilon\vert
\nabla f\vert \; u\cdot S (v).\qedhere$$
\end{proof}

According to the previous proposition, algebraic invariants of  the
linear transformation matrix in the hyperspace can be expressed by
corresponding matrix of higher order. This idea is used in \cite{Nishimura1991}
to calculate Gauss curvature of the hypersurface. Some constructions in the proof of next lemma are very similar.

\begin{lem}\label{lema algebarska}
Let $V$ be a $n-$dimensional real vector  space, $T$ a
linear transformation of hyperspace $W$ and $N$ unit vector
orthogonal to $W$. If for any two vectors $u,v\in W$ scalar
product $u\cdot T(v)$ is expressed by matrix of order $n$,
$A=[a_{ij}]$ in the following way
\begin{equation}\label{uslov lema algebarska}
u\cdot T(v)=uAv^T=a_{ij}u^i v^j,
\end{equation}
then a trace  of $T$ is given by
\begin{equation}\label{tragovi lema algebarska}
trT=trA-NAN^T.
\end{equation}
\end{lem}
\begin{proof}
Let $\{X_1,X_2,\dots,X_{n-1}\}$ be a basis of hyperspace $W$ and
$Y_{\alpha}=T(X_{\alpha}).$ Decompositions of the images $Y_{\alpha}$ over
the chosen basis
\begin{equation}\nonumber
Y_{\alpha}=b_{\beta\alpha}X_{\beta};\qquad B=\left[b_{\alpha\beta}\right]
\end{equation}
give a matrix $B$ of the linear transformation $T$ in $W$, such that
$$trT=trB.$$
Vectors $X_{\alpha},Y_{\alpha}$ actually are row matrices of dimension
$1\times n$, and $B$ is of order $n-1$. Now, one needs matrices of
order $n$:
$$\widetilde{B}=\left[\begin{array}{cc}
  B & 0 \\
  0 & 1
\end{array}\right],\qquad
 X=\left[\begin{array}{c}
  X_1 \\
  \vdots\\
  X_{n-1}\\
  N
\end{array}\right],\qquad
Y=\left[\begin{array}{c}
  Y_1 \\
  \vdots\\
  Y_{n-1}\\
  N
\end{array}\right];$$
and matrices of order $n+1$:
$$\widetilde{A}=\left[\begin{array}{cc}
  A & N^T \\
  N & 0
\end{array}\right],\qquad
\widetilde{X}=\left[\begin{array}{cc}
  X & 0 \\
  0 & 1
\end{array}\right].$$
Block matrices multiplications give:
\begin{equation}\label{XXB}
  XX^T\widetilde{B}=\left[\begin{array}{cc}
  X_{\alpha}Y_{\beta} & 0 \\
  0 & 1
\end{array}\right],
\end{equation}
\begin{equation}\label{XAX}
\widetilde{X}\widetilde{A}\widetilde{X}^T=\left[\begin{array}{ccc}
  X_{\alpha}AX^{T}_{\beta} & X_{\alpha}AN^{T} & 0 \vspace{1mm}\\
  NAX^{T}_{\beta} & NAN^T & 1 \\
  0 & 1 & 0
\end{array}\right].
\end{equation}
Without loss of generality, an assumption  of orthogonality of the basis  $\{X_1,X_2,\dots,X_{n-1}\}$ may be used. Then, $X$ and $\widetilde{X}$ are orthogonal matrices and according to properties of the trace, two previous equations  (\ref{XXB}) and (\ref{XAX}) result
\begin{equation}\nonumber
tr\left[ X_{\alpha}Y_{\beta}\right]=tr\widetilde{B}-1=trT
\end{equation}
\begin{equation}\nonumber
trA=tr\left[ X_{\alpha}AX^T_{\beta}\right]+NAN^T.
\end{equation}
These equations with (\ref{uslov lema algebarska}) give (\ref{tragovi lema algebarska}).
\end{proof}

This result is used to express a mean curvature of the hypersurface (\ref{def I}) by
its defining function.
\begin{thrm}\label{H rimanove hiperpov}
A mean curvature of the oriented hypersurface $(I,N)$ (\ref{def I}), (\ref{normal}) in Riemannian space
is
\begin{equation}\nonumber
H=\dfrac{1}{n-1}\left(tr[f_{ij}]-N[f_{ij}]N^T\right).
\end{equation}
Used notation is $f_i=\dfrac{\partial f}{\partial x^i}$, $f_{ij}=\dfrac{\partial^2 f}{\partial x^i \partial x^j}.$
\end{thrm}
\begin{proof}
Proposition \ref{propozicija} provides conditions of the lemma \ref{lema algebarska}.
By putting $a_{ij}=f_{ij}$, the proof is obvious.
\end{proof}

\section{Indicatrix of a Finsler manifold}
According to properties of fundamental function, indicatrix of the Finsler manifold $(M,F)$,
$$I_x=\left\lbrace y\in T_x M \vert F(x,y)=1\right\rbrace $$
is the convex oriented hypersurface \cite{Rund1959}\cite{Bao2000}\cite{Matsumoto1986} of Riemaniann space $T_xM$, so results of the previous section can be applied. Defining function of the hypersurface (\ref{def I}) can be presented in few forms:
\begin{eqnarray}
f(y) & =F(x,y)-1\nonumber\\
& =F^2(x,y)-1\label{def.func}\\
& =g_{ij}y^iy^j-1;\nonumber
\end{eqnarray}
and unit normal is just radius vector \cite{Kikuchi1962}\cite{Rund1959}
$$N(y)=y.$$
Second order derivative of defining function (\ref{def.func}) can be related with fundamental function or metric tensor,
\begin{equation}\nonumber
\dfrac{\partial^2 f}{\partial y^i \partial y^j}=F_{y^iy^j}=g_{ij},
\end{equation}
where an abbreviation is used, $F_{y^iy^j}=\dfrac{\partial^2 F}{\partial y^i \partial y^j}$.

Considering these facts, the mean curvature theorem  \ref{H rimanove hiperpov} can be adapted for Finsler manifold.
\begin{thrm}\label{H indikatrise}
A mean curvature of the indicatrix of Finsler manifold is
\begin{eqnarray}\nonumber
H=\dfrac{1}{n-1}\left(tr[F_{y^i y^j}]-1\right). \label{H preko F}
\end{eqnarray}
\end{thrm}

Any Riemannian space is locally Euclidean, so  in $T_x M$ a coordinates can be chosen such that
$g_{ij}=\delta_{ij}$, which means $tr[g_{ij}]=n$. A trace of the matrix is an algebraic invariant, therefore the last relation is coordinate free. In that way, a necessary condition for the fundamental function is obtained.
\begin{prop}
If $(M,F)$ is Finsler manifold, then a trace of second derivatives (Hessian) of fundamental function is just a dimension of basic manifold,
$$tr[F_{y^i y^j}]=dimM.$$
\end{prop}
Combination of the previous property with theorem \ref{H indikatrise} gives final result:
\begin{thrm}\label{H indikatrise je 1}
A mean curvature of the indicatrix of Finsler manifold is constant and identically  equal 1,  $$H=1$$.
\end{thrm}

\end{document}